\documentclass[12pt]{amsart}
\usepackage{graphicx}
\advance\hoffset by -0.5 truecm

\usepackage{amssymb}
\newtheorem{Theorem}{Theorem}[section]
\newtheorem{Lemma}[Theorem]{Lemma}

\newtheorem{Remark}[Theorem]{Remark}

\def\R{{\mathbb R}}

\def \H{{\mathbb H}}

\def \0{\lambda_{0}}

\def \om{\omega}

\def\Id{{\rm Id}}
\def\Om{\Omega}
\def\HS{\mathcal{HS}}
\def\SHS{\mathcal{SHS}}

\begin{document}
\title[]{Stability is not open}

\author[K. Cieliebak]{Kai Cieliebak}
 \address{Ludwig-Maximilians-Universit\"at, 80333 M\"unchen, Germany}
 \email {kai@math.lmu.de}

\author[U. Frauenfelder]{Urs Frauenfelder}
 \address{Department of Mathematics and Research Institute of Mathematics,
  Seoul National University}
 \email {frauenf@snu.ac.kr}

\author[G.P. Paternain]{Gabriel P. Paternain}
 \address{ Department of Pure Mathematics and Mathematical Statistics,
University of Cambridge,
Cambridge CB3 0WB, UK}
 \email {g.p.paternain@dpmms.cam.ac.uk}


\date{August 2009}

\begin{abstract} We give an example of a symplectic manifold with a
stable hypersurface such that nearby hypersurfaces are typically
unstable.
\end{abstract}

\maketitle

\section{Introduction}

A closed hypersurface $\Sigma$ in a symplectic manifold $(M,\Om)$ is
called {\em stable} if a neighbourhood of $\Sigma$ can be foliated by
hypersurfaces whose characteristic foliations are conjugate. Here the
characteristic foliation on a hypersurface $\Sigma$ is defined
by the 1-dimensional distribution $\ker(\Om|_\Sigma)$.
Stability was introduced in~\cite{HZ} as a
condition on hypersurfaces for which the Weinstein conjecture can be
proved. More recently, it has attained importance as the condition
needed for the compactness results underlying Symplectic Field
Theory~\cite{EGH,BEHWZ,CM} and Rabinowitz Floer
homology~\cite{CF,CFP}.
\smallskip

Let us consider, in a fixed symplectic manifold $(M,\Om)$, the space
$\HS$ of closed hypersurfaces equipped with the $C^\infty$-topology
and its subset $\SHS$ of stable hypersurfaces. It is easy to see that
$\SHS$ is not closed: For example, the horocycle flow on a hyperbolic
surface defines a hypersurface which is unstable but the smooth limit
of stable ones; see~\cite{CFP} for many more examples. On the other
hand, $\SHS$ contains open components, e.g.~those corresponding to
hypersurfaces of contact type. This prompted the question whether the
set $\SHS$ is actually open in $\HS$. The result of this paper shows
that this is not the case.
\smallskip

\begin{Theorem}\label{thm:main}
There exists a stable closed hypersurface $\Sigma$ in a symplectic
6-manifold such that nearby hypersurfaces are typically unstable in
the following sense: There exists a neighbourhood of $\Sigma$ in $\HS$
which contains an open dense set consisting of unstable hypersurfaces.
\end{Theorem}
\smallskip

The theorem continues to hold if the $C^\infty$ topology is replaced
by the $C^k$ topology for some $k\geq 2$ and hypersurfaces are only
assumed to be of class $C^k$.
\smallskip

The theorem can be rephrased in terms of {\em stable Hamiltonian
  structures}~\cite{BEHWZ,CM,CV}.
A two-form $\omega$ on an odd-dimensional manifold $\Sigma$ is
called a {\em Hamiltonian structure} if it is closed and maximally
nondegenerate in the sense that its kernel distribution is
one-dimensional. It is called {\em stable} if there exists a one-form
$\lambda$ such that $\lambda|_{\ker\om}\neq 0$ and $\ker \omega
\subset \ker d\lambda$. Then a hypersurface $\Sigma$ in a symplectic
manifold $(M,\Om)$ is stable iff $\Om|_\Sigma$ defines a stable
Hamiltonian structure, and every stable Hamiltonian structure arises
as a stable hypersurface in some symplectic manifold~\cite{CM}. Now
Theorem~\ref{thm:main} can be rephrased as follows: {\sl There exists
a stable Hamiltonian structure $\om$ on a closed 5-manifold $\Sigma$
such that nearby Hamiltonian structures with the same cohomology class as
$\om$ are typically unstable.}
\smallskip

Theorem~\ref{thm:main} has implications on the foundations of
holomorphic curve theories such as Symplectic Field
Theory~\cite{EGH,BEHWZ,CM} and Rabinowitz Floer
homology~\cite{CF,CFP}.
For the construction of those theories one needs to perturb a given
stable Hamiltonian structure to make all closed characteristics
nondegenerate. Theorem~\ref{thm:main} suggests that such a
perturbation may not be possible within the class of stable
Hamiltonian structures (see also~\cite{CV} for a result pointing in
the same direction).
In Rabinowitz Floer homology this problem can be overcome in the
following way~\cite{CFP}:  One chooses an additional Hamiltonian
perturbation of the Rabinowitz action functional. For a generic
small perturbation the Rabinowitz action functional becomes Morse, but
for the perturbed action functional one might lose compactness. However,
one can still define a boundary operator by taking into account
only gradient flow lines close to the original ones. We
wonder if a similar strategy can be applied to SFT as well.

\section{Preliminaries on Anosov Hamiltonian structures}

\noindent{\bf Anosov Hamiltonian structures. }
Recall that the flow $\phi_t$ of a vector field $F$ on a closed manifold
$\Sigma$ is {\em Anosov} if there is a splitting $T\Sigma=\R
F\oplus E^s\oplus E^u$ and positive constants $\lambda$ and $C$ such that for
all $x\in \Sigma$
\[|d_{x}\phi_{t}(v)|\leq
C e^{-\lambda t}|v|\;\;\mbox{\rm for}\;v\in E^s\;\;\mbox{\rm and}\;t\geq 0,\]
\[|d_{x}\phi_{-t}(v)|\leq
C e^{-\lambda t}|v|\;\;\mbox{\rm for}\;v\in E^u\;\;\mbox{\rm and}\;t\geq 0.\]
If an Anosov vector field $F$ is rescaled by a positive function its
flow remains Anosov~\cite{AS,Pa}. It will be useful
for us to know how the bundles $E^s$ and $E^u$ change when
we rescale $F$ by a smooth positive function $r:\Sigma\to\R_{+}$.
Let $\tilde\phi$ be the flow of $rF$ and $\tilde{E}^s$ its stable bundle.
Then (cf. \cite{Pa})
\begin{equation}
\tilde{E}^{s}(x)=\{v+z(x,v)F(x):\;\;\;v\in E^{s}(x)\},
\label{eq:stc}
\end{equation}
where $z(x,v)$ is a continuous 1-form (i.e. linear in $v$ and continuous
in $x$). Moreover, if we let $l=l(t,x)$ be (for fixed $x$) the inverse
of the diffeomorphism
\[t\mapsto\int_{0}^{t}r(\phi_{s}(x))^{-1}\,ds\]
then
\begin{equation}
d\tilde\phi_{t}(v+z(x,v)F(x))=
d\phi_{l}(v)+z(\phi_{l}(v),d\phi_{l}(v))F(\phi_{l}(x)).
\label{eq:tch}
\end{equation}
This shows that for closed $\Sigma$ the flow $\tilde\phi_t$ is again
Anosov. There is a similar expression for $\tilde{E}^{u}$.
It is clear from the discussion above that the weak bundles $\R F\oplus E^s$
 and
$\R F\oplus E^u$ do not change under rescaling of $F$ (the
strong bundles $E^{s,u}$ are indeed affected by rescaling as we have just seen).

Let $(\Sigma,\omega)$ be a Hamiltonian structure. We say that the structure is
{\em Anosov} if the flow of any vector field $F$ spanning $\ker\om$ is Anosov.

We say that an Anosov Hamiltonian structure satisfies the {\it
  $1/2$-pinching condition} or that it is {\it 1-bunched} \cite{H,H2} if
for any vector field $F$ spanning $\ker\om$ with flow $\phi_t$
there are functions $\mu_{f}, \mu_{s} :\Sigma\times \R_{+}\to \R_{+}$
such that
\begin{itemize}
\item $\lim_{t\to\infty}\sup_{x\in\Sigma}\frac{\mu_{s}(x,t)^{2}}{\mu_{f}(x,t)}=0$;
\item $\mu_{f}(x,t)|v|\leq |d\phi_{t}(v)|\leq \mu_{s}(x,t)|v|$ for all
$x\in\Sigma$, $t>0$ and $v\in E^{s}(x)$, and $\mu_{f}(x,t)|v|\leq |d\phi_{-t}(v)|\leq \mu_{s}(x,t)|v|$ for all
$x\in\Sigma$, $t>0$ and $v\in E^{u}(\phi_{t}x)$.
\end{itemize}
We remark that the 1/2-pinching condition is invariant under rescaling.
Indeed, consider the flow $\tilde\phi_{t}$ of $rF$. It is clear from
(\ref{eq:stc}) and (\ref{eq:tch}) that there is a positive constant $\kappa$
such that
\[\frac{1}{\kappa}\mu_{f}(x,l(t,x))|\tilde v|\leq
|d\tilde\phi_{t}(\tilde v)|\leq \kappa \mu_{s}(x,l(t,x))|\tilde v|\]
for $t>0$ and $\tilde v\in \tilde{E}^{s}$ (with a similar expression
for $\tilde{E}^{u}$). We know that given $\varepsilon>0$, there exists
$T>0$ such that
for all $x\in\Sigma$ and all $t>T$ we have
\[\frac{\mu_{s}(x,t)^{2}}{\mu_{f}(x,t)}<\varepsilon.\]
On the other hand, there exists $a>0$ such that $l(t,x)\geq at$ for
all $x\in\Sigma$ and $t>0$. Hence, if we choose $t>T/a$ we have
\[\frac{\mu_{s}(x,l(t,x))^{2}}{\mu_{f}(x,l(t,x))}<\varepsilon\]
for all $x\in\Sigma$.
Therefore
$$\lim_{t\to\infty}\sup_{x\in\Sigma}\frac{\mu_{s}(x,l(t,x))^{2}}{\mu_{f}(x,l(t,x))}=0$$
and thus $\tilde\phi_{t}$ is also $1/2$-pinched.

Hence the Anosov
property as well as the $1/2$-pinching condition are invariant under
rescaling and thus intrinsic properties of the Hamiltonian structure.
One of the main consequences of the $1/2$-pinching condition is that
the weak bundles $\R F\oplus E^s$ and $\R F\oplus E^u$ are of class
$C^1$ \cite[Theorem 5]{H2} (see also \cite{HPS}).
\medskip

\noindent{\bf Stable Anosov Hamiltonian structures. }
Suppose now $(\Sigma,\om)$ is a {\it stable} Anosov Hamiltonian structure
satisfying the $1/2$-pinching condition. Let $\lambda$ be a
stabilizing 1-form and $R$ the Reeb vector field defined by
$i_R\om=\0$ and $\lambda(R)=1$. Invariance under the
flow implies that $\om$ and $\lambda$ both vanish on $E^s$ and $E^u$.
Since the flow $\phi_t$ of $R$ is Anosov and
$E^s\oplus E^u=\ker\lambda$ which is $C^\infty$, it follows that
$E^s=\ker\lambda\cap(\R F\oplus E^s)$ and $E^u$ must be $C^1$. Under
these conditions we can introduce
the {\it Kanai connection} \cite{Ka} which is defined as follows.

Let $I$ be the $(1,1)$-tensor on $\Sigma$ given by $I(v)=-v$ for $v\in E^s$,
$I(v)=v$ for $v\in E^u$ and $I(R)=0$. Consider the symmetric non-degenerate
 bilinear form given by
\[h(X,Y):=\om(X,IY)+\lambda\otimes\lambda(X,Y).\]
The pseudo-Riemannian metric $h$ is of class $C^1$ and thus there exists
a unique $C^0$ affine connection $\nabla$ such that:
\begin{enumerate}
\item $h$ is parallel with respect to $\nabla$;
\item $\nabla$ has torsion $\om\otimes R$.
\end{enumerate}
This connection has the following desirable properties \cite{Fe,Ka}:
it is invariant under $\phi_t$
and the Anosov splitting is invariant under $\nabla$: if $X$ is any section
of $E^{s,u}$, $\nabla_{v}X\in E^{s,u}$ for any $v$.

The other good consequence of the $1/2$-pinching condition, besides
$C^1$ smoothness of the bundles, is the following lemma (cf. \cite[Lemma 3.2]{Ka}).

\begin{Lemma} $\nabla(d\lambda)=0$.
\label{parallel}
\end{Lemma}

\begin{proof}
Suppose $\tau$ is any invariant $(0,3)$-tensor annihilated by $R$.
We claim that $\tau$ must vanish. To see this, consider for example
a triple of vectors $(v_1,v_2,v_3)$ where $v_1,v_2\in E^{s}$
but $v_{3}\in E^{u}$. Then there is a constant $C>0$ such that for all
$t\geq 0$
\begin{align*}
|\tau_{x}(v_{1},v_{2},v_{3})|&=|\tau_{\phi_{t}x}(d\phi_{t}(v_{1}),d\phi_{t}(v_{2}),d\phi_{t}(v_{3}))|\\
&\leq C {\mu_{s}(x,t)}^{2}\mu_{f}(x,t)^{-1}|v_{1}||v_{2}||v_{3}|.
\end{align*}
By the $1/2$-pinching condition the last expression tends to zero
as $t\to\infty$ and therefore $\tau_{x}(v_{1},v_{2},v_{3})=0$.
The same will happen for other possible triples $(v_1,v_2,v_3)$
when we let $t\to\pm\infty$.

Since $d\lambda$ and $\nabla$ are $\phi_t$-invariant, so is
 $\nabla(d\lambda)$. Since $i_{R}d\lambda=0$, $\nabla(d\lambda)$ is
also annihilated by $R$ (to see that $\nabla_{R}(d\lambda)=0$ use
that $d\lambda$ is $\phi_t$-invariant and that $\nabla_{R}=L_{R}$).
Hence by the previous argument applied to $\tau=\nabla(d\lambda)$
we conclude that $\nabla(d\lambda)=0$ as desired.
\end{proof}

\medskip

\noindent{\bf Quasi-conformal Anosov Hamiltonian structures. }
Let $\phi_t$ be an Anosov flow on $\Sigma$ endowed with a $C^0$-Riemannian metric. Consider the following functions on $\Sigma\times\R$:
\[K^s(x,t)=\frac{\max\{|d\phi_{t}(v)|:\;v\in E^s(x),\;|v|=1\}}{\min\{|d\phi_{t}(v)|:\;v\in E^s(x),\;|v|=1\}},\]
\[K^u(x,t)=\frac{\max\{|d\phi_{t}(v)|:\;v\in E^u(x),\;|v|=1\}}{\min\{|d\phi_{t}(v)|:\;v\in E^u(x),\;|v|=1\}}.\]
The flow $\phi_t$ is said to be {\it quasi-conformal} if $K^u$ and $K^s$ are both bounded on
$\Sigma\times \R$. This property is clearly independent of the choice of Riemannian metric
used to define $K^s$ and $K^u$. Moreover it is shown in \cite[Proposition 3.5]{S}
that quasi-conformality is independent of times changes, thus it makes sense to talk
about quasi-conformal Anosov Hamiltonian structures. The next theorem will be useful for us.

\begin{Theorem}[\cite{S}, Theorems 1.3 and 1.4]
Let $\phi_t$ be a topologically mixing Anosov flow with $\mbox{\rm
  dim}\,E^s\geq 2$ and $\mbox{\rm dim}\,E^u\geq 2$.
If $\phi_t$ is quasi-conformal, then the weak bundles are $C^{\infty}$.
\label{sadovs}
\end{Theorem}

Recall that $\phi_t$ is topologically mixing if for any two nonempty open sets $U$ and $V$ in $\Sigma$, there
is a compact set $K\subset \mathbb R$ such that for every $t\in \mathbb R\setminus K$ we have
$\phi_{t}(U)\cap V\neq\emptyset$. Recall also that $\phi_t$ is said to be transitive if there
is a dense orbit. Our Anosov flows will always be transitive since they preserve a smooth volume form
\cite[Chapter 18]{KH}.

\section{A theorem}

\begin{Theorem} Let $(\Sigma, \omega)$ be a 1/2-pinched Anosov Hamiltonian structure
with $[\omega]\neq 0$, but $[\om^2]=0$. Suppose in addition that $\Sigma$ fibres over a closed 3-manifold
with fibres diffeomorphic to $S^2$ and transversal to the weak subbundles.
Then, if $(\Sigma, \omega)$ is stable,
the weak subbundles must be $C^{\infty}$.
\label{theorem:key}
\end{Theorem}

\begin{proof} The proof of this theorem is very much inspired by the proof
of Theorem 2 in \cite{Ka}.
We first make the following observation:
\begin{itemize}
\item $E^s$ ($E^u$) cannot contain a nontrivial proper continuous subbundle.
\end{itemize}
Indeed since $\R R\oplus E^u$ is transversal to the fibres of the fibration $\Sigma\to M$
by 2-spheres, we can write $T\Sigma=V\oplus \R R\oplus E^u$ where
$V$ is the vertical subbundle of the fibration. Using this splitting we may define
an isomorphism $E^s\mapsto V$ and since the tangent bundle of $S^2$ does not admit
a nontrivial proper continuous subbundle, the same holds for $E^s$ (and $E^u$).

Next we observe that the stabilizing 1-form $\lambda$ cannot be
closed. Indeed, write $\om^2=d\tau$ and note that
if $\lambda$ was closed, then the volume form $\lambda\wedge d\tau$
would be exact, which is absurd.

Since $\om$ is non-degenerate, there exists a smooth bundle map
$L:E^s\oplus E^u\to E^s\oplus E^u$ such that for sections $X,Y$
of $E^s\oplus E^u$
\[d\lambda(X,Y)=\om(LX,Y)=\om(X,LY).\]
The map $L$ is invariant under $\phi_t$ and preserves the
 decomposition $E^s\oplus E^u$, i.e.
 $L=L^s+L^u$, where $L^s:E^s\to E^s$ and $L^u:E^u\to E^u$. In
 particular, $L$ commutes with $I$.
By Lemma \ref{parallel}, the $1/2$-pinching condition implies that $\nabla(d\lambda)=0$ and thus
$L$ is parallel with respect to $\nabla$.
Note that by transitivity
of $\phi_t$, the characteristic polynomial of $L_{x}^s$
is independent of $x\in\Sigma$. Let $\rho\in\mathbb C$ be an eigenvalue
of $L^s$. Consider $A:=L^s-\Re(\rho)\Id$. Note that $A$ cannot be zero: Otherwise
$d\lambda=c\,\omega$ for a constant $c\in\R$; since $\lambda$ is not closed, $c\neq 0$, which in turns implies $[\omega]=0$, contradicting the hypotheses of the theorem.

Clearly $A^2$ has $\mu:=-\Im(\rho)^2$ as an eigenvalue.
Let $H\subset E^s$ denote the eigenspace of the eigenvalue $\mu$. Since $L^s$ is parallel
it has the same dimension at every point $x\in \Sigma$ and
since $E^s$ cannot contain a nontrivial proper continuous subbundle, we deduce that
$H=E^s$. Hence $A^2=\mu \Id$. Moreover $\mu\neq 0$, otherwise $\mbox{\rm ker} A$ would
be a nontrivial proper continuous subbundle of $E^s$.
Therefore we have proved that
\[\mathbb J^s:=\frac{1}{\Im(\rho)}(L^s-\Re(\rho)\Id),\]
defines a parallel almost complex structure on $E^s$ of class $C^1$ invariant
under $\phi_t$. Similarly we obtain an almost complex structure
$\mathbb J^u$ on $E^u$.

Now choose a Riemannian metric on $E^s$ (resp. $E^u$) which is invariant under $\mathbb J^s$
(resp. $\mathbb J^u$). By declaring $E^s$, $E^u$ and $\R R$ orthogonal and $R$ with norm 1, we obtain a metric (of class $C^1$) on $\Sigma$ such that with respect to this metric
\[\frac{\max\{|d\phi_{t}(v)|:\;v\in E^s(x),\;|v|=1\}}{\min\{|d\phi_{t}(v)|:\;v\in E^s(x),\;|v|=1\}}=1,\]
for all $t\in\R$ and $x\in \Sigma$. This is because $\phi_{t}$
preserves $\mathbb J^s$ and $E^s$ has rank two.
Similarly for $E^u$. This shows that $(\Sigma,\omega)$ is a quasi-conformal Anosov Hamiltonian
structure.

Finally we note that if a transitive Anosov flow is not topologically mixing, then
by a theorem of J. Plante \cite{Pl} it must be a suspension with constant return function. In particular, this
implies that there is
a closed 1-form $\beta$ such that $\beta(R)>0$. The same argument above that proved that
$\lambda$ cannot be closed shows that such a $\beta$ cannot exist. Hence $\phi_t$
is topologically mixing and by Theorem \ref{sadovs} the weak bundles must be $C^{\infty}$.
\end{proof}

\begin{Remark}{\rm
Note that the proof above only requires $\lambda$ to be of class $C^2$.
}
\end{Remark}

\section{The example} Let $\Gamma$ be a discrete group of isometries of $\H^3$ such that
$M:=\Gamma\setminus\H^3$ is a closed orientable hyperbolic 3-manifold.
We consider the geodesic flow acting on the unit sphere bundle $SM$ and let $\alpha$
be the canonical contact 1-form.

The space of invariant 2-forms of the geodesic flow of $M=\Gamma\setminus \H^3$
has dimension two~\cite[Claim 3.3]{Ka}.
It is spanned by the 2-form $d\alpha$, where $\alpha$ is the
canonical contact form on the unit sphere bundle $SM$, and the following additional
2-form $\psi$ which we now describe.
Given a unit vector $v\in T_{x}\H^3$, let $i(v):T_{x}\H^3\to T_{x}\H^3$ be the linear map
defined by $i(v)(v)=0$ and $i(v)$ rotates vectors in $\{v\}^{\perp}$ by $\pi/2$ according
to the orientation of $\H^3$. Any vector $\xi\in T_{v}S\H^3$ can be written
as $\xi=(\xi_{H},\xi_{V})$ with the usual identification of horizontal and vertical components (cf. \cite{Pat}).
Define $J_{v}:T_{v}S\H^3\to T_{v}S\H^3$ as
\begin{equation}
J_{v}(\xi_{H},\xi_{V})=(i(v)\xi_{V},i(v)\xi_{H}).
\label{eq:J}
\end{equation}
Then
\begin{equation}
\psi_{v}(\xi,\eta):=d\alpha_{v}(J_{v}\xi,\eta)=\langle i(v)\xi_{V},\eta_{V}\rangle-\langle i(v)\xi_{H},\eta_{H}\rangle.
\label{eq:psi}
\end{equation}
Clearly this construction descends to $SM$ where we use the same notation
($\psi$, $\alpha$, etc.)
In a moment we will check that $\psi$ is invariant under $\phi_t$, but before we do so, let us describe the stable and unstable bundles of $\phi_t$ and the action of $d\phi_t$ on them.
Recall that $d\phi_{t}(\xi_{H},\xi_{V})=(Y(t),\dot{Y}(t))$
where $Y$ is the unique Jacobi field (along the geodesic $\pi\phi_t(v)$, where
$\pi:SM\to M$ is foot-point projection) with initial conditions $(\xi_{H},\xi_{V})$.
Solving the Jacobi equation $\ddot{Y}-Y=0$ we find:
\[E^s(v)=\{(w,-w):\;w\perp v\},\]
\[E^u(v)=\{(w,w):\;w\perp v\}.\]
Note that $J$ leaves $E^s$ and $E^u$ invariant. Moreover
\[d\phi_{t}(w,-w)=e^{-t}(e_{w}(t),-e_{w}(t)),\]
\[d\phi_{t}(w,w)=e^{t}(e_{w}(t),e_{w}(t)),\]
where $e_{w}(t)$ is the parallel transport of $w$ along the geodesic $\pi\phi_t(v)$.
Since $e_{i(v)w}(t)=i(\pi\phi_{t}v)e_{w}(t)$ we see that $d\phi_{t}$ preserves $J$.
Since $d\alpha$ is also $\phi_t$ invariant, it follows that $\psi$ is invariant.
Note that $i_R\psi=0$ for the Reeb vector field $R$ of $\alpha$.

\begin{Lemma} The invariant 2-form $\psi$ is closed but not exact.
The 4-form $\psi^2$ is exact and $(SM,\psi)$ is a stable Hamiltonian
structure
with stabilizing 1-form $\alpha$ and Reeb vector field $R$.
\label{lemma:psi}
\end{Lemma}

\begin{proof} The 3-form $d\psi$ is invariant under $\phi_t$ and is annihilated by $R$.
Then the proof of Lemma \ref{parallel} shows that $d\psi=0$ (obviously $\phi_t$ is 1/2-pinched).
In order to show that $[\psi]\neq 0$, consider $S_x$ the 2-sphere of unit vectors
in $T_{x}\H^3$. A tangent vector $\xi\in T_{v}S_{x}$ has the form $\xi=(0,w)$ where
$w\perp v$. If we take two tangent vectors $\xi=(0,w)$, $\eta=(0,u)\in T_{v}S_{x}$,
from (\ref{eq:J}) and (\ref{eq:psi}) we see that
\[\psi_{v}(\xi,\eta)=\langle i(v)w,u\rangle.\]
This implies that
\[\int_{S_{x}}\psi\neq 0\]
and thus $[\psi]\neq 0$.
Consider now the invariant 4-form $\psi^2$ and the invariant 5-form $\alpha\wedge\psi^2$.
By transitivity, there is a constant $k$ such that $\alpha\wedge\psi^2=k\,\alpha\wedge(d\alpha)^2$.
Contracting with $R$ we see that $\psi^2$ must be
$k\,(d\alpha)^2$ and therefore exact.
Finally, it is immediate from the definition (\ref{eq:psi}) of $\psi$ that its restriction to
$E^s\oplus E^u=\mbox{\rm ker}\,\alpha$ is non-degenerate. Hence $(SM,\psi)$ is a Hamiltonian
structure with stabilizing 1-form $\alpha$ and Reeb vector field $R$.
\end{proof}

Now let $X:=SM\times (-\varepsilon,\varepsilon)$ and $\tau:X\to SM$ the obvious projection.
Define $\omega_{X}:=d(r\tau^*\alpha)+\tau^*\psi$, where $r\in (-\varepsilon,\varepsilon)$.
For $\varepsilon$ small enough $(X,\omega_{X})$ is a symplectic manifold and
$r=0$ is the stable hypersurface $(SM,\psi)$.

We have now come to our main result which implies
Theorem~\ref{thm:main} in the introduction.

\begin{Theorem} A typical hypersurface $\Sigma\subset X$ near $SM$ is not stable.
\end{Theorem}

\begin{proof} Consider a hypersurface $\Sigma$ near $r=0$ and let $\omega$ be $\omega_{X}$
restricted to $\Sigma$. By Lemma \ref{lemma:psi}, $[\omega]\neq 0$, but $[\omega^2]=0$.
Since $SM$ fibres over $M$ with fibres given by 2-spheres transveral to the weak bundles
the same holds true for $\Sigma$ (recall that under perturbations the stable and unstable
bundles vary continuously).
Finally we note that $(\Sigma,\omega)$ is 1/2-pinched. Indeed, recall that for the geodesic flow
of $M$, we have
\[|d\phi_{t}(\xi)|=e^{-t}|\xi|\;\mbox{\rm for}\;\xi\in E^{s},\]
\[|d\phi_{t}(\xi)|=e^{t}|\xi|\;\mbox{\rm for}\;\xi\in E^{u}.\]
Thus for a flow $\varphi_t$ which is $C^1$ close to $\phi_t$ we get
\[\frac{1}{C}|\xi|e^{-At}\leq|d\varphi_{t}(\xi)|\leq C|\xi|e^{-at}\;\;\mbox{\rm for}\;\xi\in E^s\;\;\mbox{\rm and}\;t\geq 0,\]
\[\frac{1}{C}|\xi|e^{-At}\leq|d\varphi_{-t}(\xi)|\leq C|\xi|e^{-at}\;\;\mbox{\rm for}\;\xi\in E^u\;\;\mbox{\rm and}\;t\geq 0,\]
where all the constants $C,A,a$ are close to 1. Thus $(\Sigma,\omega)$ is 1/2-pinched.

We can now apply Theorem \ref{theorem:key} to conclude that if $\Sigma$ near $r=0$ is stable, then
the weak bundles must be $C^\infty$.
However, a theorem of Hasselblatt \cite[Corollary 1.10]{H} asserts that an open and dense set
of symplectic Anosov systems does not have weak bundles of class
$C^{2-\varepsilon}$. Thus a typical hypersurface $\Sigma$ near $r=0$
cannot be stable.
\end{proof}

\begin{Remark}{\rm It is possible to prove the last theorem without appealing to
Theorem \ref{sadovs}. An inspection of the proof of Theorem \ref{theorem:key}
shows that since $d\phi_t$ preserves $\mathbb J$, all the closed orbits
are actually 2-bunched in the terminology of \cite{H}, and the local perturbation argument
in \cite[Section 4]{H} implies that an open and dense set of symplectic Anosov systems
 does not have all closed orbits being $2$-bunched (this fact is actually used in
 the proof of \cite[Corollary 1.10]{H} quoted above). 
Of course, the conclusion of Theorem \ref{theorem:key} is stronger if we use
 Theorem \ref{sadovs}.

}
\end{Remark}

\end{document}